\documentclass[12pt]{amsart}

\usepackage{melanie}

\usepackage{fullpage}

\newtheorem{nts}{Note to self}

\title{The distribution of the number of points on trigonal curves over $\F_q$}

\author{Melanie Matchett Wood}
\address{Department of Mathematics\\
University of Wisconsin \\ 480 Lincoln Drive \\
Madison, WI 53705 USA\\
and
American Institute of Mathematics\\360 Portage Ave \\
Palo Alto, CA 94306-2244 USA} 
\email{mmwood@math.wisc.edu}

\begin{document}
\maketitle

\begin{abstract}
We give a short determination of the distribution of the number of $\F_q$-rational points on a random trigonal curve over $\F_q$, in the limit as the genus of the curve goes to infinity.
In particular, the expected number of points is $q+2-\frac{1}{q^2+q+1}$, contrasting with recent analogous results for cyclic $p$-fold covers of $\P^1$ and plane curves
 which have an expected number of points of $q+1$ (by work of Kurlberg, Rudnick, Bucur, David, Feigon and Lal\'in) and curves which are complete intersections
which have an expected number of points $<q+1$ (by work of Bucur and Kedlaya).  We also give a conjecture for the expected number of points on a random $n$-gonal curve with full $S_n$ monodromy
based on function field analogs of Bhargava's heuristics for counting number fields.
\end{abstract}

\section{Introduction}
If we fix a finite field $\F_q$, we can ask about the distribution of the number of ($\F_q$-rational) points on a random curve over $\F_q$.  There has been a surge of recent activity on this question including definitive answers of Kurlberg and Rudnick  for hyperelliptic curves \cite{KuRu}, of Bucur, David, Feigon and Lal\'in for cyclic $p$-fold covers of $\P^1$ \cite{BuDaFeLa1,BuDaFeLa2}, of Bucur, David, Feigon and Lal\'in for plane curves \cite{BuDaFeLa3}, and of Bucur and Kedlaya on curves that are complete intersections in smooth quasiprojective subschemes of $\P^n$ \cite{BuKe}.  In the first three cases, the average number of points on a curve in the family is $q+1$.  In contrast, for curves that are complete intersections in $\P^n$, the average number of points is $<q+1$, despite, as pointed out by Bucur and Kedlaya, the abundance of $\F_q$ points lying around in $\P^n$.  In this paper, we give the distribution of the number of points on trigonal curves over $\F_q$ (i.e. curves with a degree $3$ map to $\P^1$), and in particular show that the average number of points is greater that $q+1$.

Let 
$$T_g:=\{\pi:C\ra \P^1 | C \textrm{ is a smooth, geometrically integral, genus $g$ curve with $\pi$ degree 3}\}.$$
Our main theorem is the following.
\begin{theorem}\label{T:Main}
Let $\F_q$ have characteristic $\geq 5$.
We have
$$\lim_{g\ra\infty} \frac{\#|\{C\in T_g(\F_q) \mid \#|C(\F_q)|=k  \}|}{\#|T_g(\F_q)|}=
\operatorname{Prob}(X_1+\dots+X_{q+1}=k),
$$
where the $X_i$ are independent identically distributed random variables and
$$
X_i=
\begin{cases}
0 &\textrm{ with probability } \frac{2q^2}{6q^2+6q+6} \\
1 &\textrm{ with probability } \frac{3q^2+6}{6q^2+6q+6}\\
2 &\textrm{ with probability } \frac{6q}{6q^2+6q+6}\\
3 &\textrm{ with probability } \frac{q^2}{6q^2+6q+6}.
\end{cases}
$$
\end{theorem}
Moreover, we give a rigorous explanation of the random variables $X_i$. 
 Theorem~\ref{T:Main} is a corollary of the following, which gives the distribution of the
 number of points in the fiber over a given $\F_q$-rational point of $\P^1$.  
\begin{theorem}\label{T:ByPoint}
Let $\F_q$ have characteristic $\geq 5$.
Given a point $z\in\P^1(\F_q)$,
$$\lim_{g\ra\infty} \frac{\#|\{(C,\pi)\in T_g(\F_q) \mid \#|\pi^{-1}(x)(\F_q)|=k  \}|}{\#|T_g(\F_q)|}=
\begin{cases}
\frac{2q^2}{6q^2+6q+6} &\textrm{ for } k=0\\
\frac{3q^2+6}{6q^2+6q+6} &\textrm{ for } k=1 \\
\frac{6q}{6q^2+6q+6} &\textrm{ for } k=2 \\
\frac{q^2}{6q^2+6q+6} &\textrm{ for } k=3 .
\end{cases}
$$
Moreover, these probabilities (of various size fibers) are independent at the 
$\F_q$ points $z_1,\dots,z_{q+1}$ of $\P^1$.
\end{theorem}

\begin{corollary}
The average number of points of a random trigonal curve over $\F_q$ (in the $g\ra \infty$ limit as above) is $q+2-\frac{1}{q^2+q+1}$.
\end{corollary}

The method used in this paper will be to relate trigonal curves to cubic extensions of function fields, and then to use the work of Datskovsky and Wright \cite{DaWr} to count cubic extensions with every possible fiberwise behavior above each rational point of the base curve.
In fact, our methods work with any smooth curve $E$ over $\F_q$ replacing $\P^1$.  
Let $T_{E,g}$ be the moduli space of genus $g$ curves with a specified degree $3$ map to $E$.
\begin{theorem}\label{T:EByPoint}
Let $\F_q$ have characteristic $\geq 5$.
Given a point $z\in E(\F_q)$,
$$\lim_{g\ra\infty} \frac{\#|\{(C,\pi)\in T_{E,g}(\F_q) \mid \#|\pi^{-1}(x)(\F_q)|=k  \}|}{\#|T_g(\F_q)|}=
\begin{cases}
\frac{2q^2}{6q^2+6q+6} &\textrm{ for } k=0\\
\frac{3q^2+6}{6q^2+6q+6} &\textrm{ for } k=1 \\
\frac{6q}{6q^2+6q+6} &\textrm{ for } k=2 \\
\frac{q^2}{6q^2+6q+6} &\textrm{ for } k=3 .
\end{cases}
$$
Moreover, these probabilities (of various size fibers) are independent at the 
$\F_q$ points of $E$.  In particular, the average number of points of a 
random curve over $\F_q$ with a degree $3$ map to $E$ (in the $g\ra \infty$ limit) is $\#|E(\F_q)|(1+\frac{q}{q^2+q+1})$.
\end{theorem}

\subsection{Related Work}
Studying the distribution of the number of points on a curve is equivalent to studying the distribution of the trace of Frobenius on the $\ell$-adic cohomology group $H^1$ .  
Bucur, David, Feigon and Lal\'in, for cyclic $p$-fold covers of $\P^1$, in fact give the finer information of the distribution of the trace of Frobenius on each subspace of $H^1$ invariant under the cyclic action.  Their methods, under appropriate interpretation, use Kummer theory to enumerate the cyclic $p$-fold covers of $\P^1$, and thus require the hypothesis that $q \equiv 1 \pmod{p}$.  The hardest part of the method is sieving for $p$-power free polynomials.

The work  of Bucur, David, Feigon and Lal\'in for plane curves \cite{BuDaFeLa3} and of Bucur and Kedlaya on curves that are complete intersections \cite{BuKe} is also given ``fiberwise,'' (as in Theorems~\ref{T:ByPoint} and \ref{T:EByPoint} of this paper) in that it computes the probability that any point is the ambient space is in a random smooth curve in the specified family and shows these probabilities are independent.  (Here we think of the embedding of a curve $i:C\ra X$, and \cite{BuDaFeLa3,BuKe} computes
for each point $z\in X(\F_q)$ the distribution of the size of the fiber $i^{-1}(z)$ for random $C$.  Note that since $i$ is an embedding, in this case the fiber has either 0 or 1 points, so these are Bernoulli random variables.)

Kurlberg and Wigman \cite{KuWi} have also studied the distribution of the number of $\F_q$ points in families of curves in which the average number of points is infinite.

In Section~\ref{S:Further} of this paper, we discuss distribution of points on other families of $n$-gonal curves, and in particular give a conjecture for the expected number of points on a random $n$-gonal curve with full $S_n$ monodromy, based on function field analogs of Bhargava's heuristics for counting number fields.  The conjecture in fact predicts the expected number of points in a fixed fiber over $\P^1$.

When one instead fixes a genus $g$ and lets $q$ tend to infinity, the philosophy of Katz and Sarnark \cite{ks} predicts how the average number of points behaves, and in particular it should be governed by statistics of random matrices in a group depending on the monodromy of the moduli space of curves under consideration.

\section{Proof of Theorem~\ref{T:EByPoint}}

Let $\F_q$ have characteristic $>3$.  Let $E$ be a smooth curve over $\F_q$, and let $k$ be the function field of $E$.  
Then smooth, integral curves over $\F_q$ with finite, degree 3 maps to $\P^1_{\F_q}$ are in one-to-one correspondence with cubic extensions $k'$ of $k$.
However, these smooth, integral curves may not be geometrically integral.  In this case, the only such possibility is given by a cubic extension of the constant field of $E$, and since we will take a $g\ra \infty$ limit, we can ignore this extension altogether.

We now consider all the possible completions of a cubic extension $k'$ of $k$.  For a place $v$ of $k$, the absolute tame Galois group of the local field $k_v$ is topologically generated by $x,y$ with the relation
$$
xyx^{-1}=y^q,
$$
where $y$ is a generator of the inertia subgroup and $x$ is a Frobenius element.  
For a place $v$ of $k$, the following is a chart of all possible cubic \'{e}tale $k_v$ algebras $L$ (and thus of all possible isomorphism types for $k'\tensor_k k_v$).
To each isomorphism type, we note the images of $x,y$ from the absolute Galois group of $k_v$ in the associated representation to $S_3$ (given by the action of Galois on the three homomorphisms $L\ra \bar{k_v}$), the number of $\F_q$ rational points in the fiber above $v$, and a constant $c_L$ that will be important later.  (It turns out that $c_L=\frac{1}{|\Aut(L)||D_{L/k_v}|})$, but this fact will not be used.)

\begin{tabular}{l|l|l|l|l}
 &  & & \# of $\F_q$   & \\[2pt]
$L$ & $x$ & $y$  &  points & $c_L$ \\[2pt]
 &  &  &  in fiber & \\[2pt]
\hline
\hline
 & & & &\\[-8pt]
$k_v^{\oplus 3}$ & $()$ & $()$ & 3 & $1/6$ \\[4pt]
\hline
 & & & &\\[-8pt]
$K\oplus k_v$, $K/k_v$ degree $2$ unram. field extn. & $(12)$ & $()$ & 3 & $1/2$ \\[4pt]
\hline
 & & & &\\[-8pt]
$K$, $K/k_v$ degree $3$ unram. field extn. & $(123)$ & $()$ & 0 & $1/3$ \\[4pt]
\hline
 & & & &\\[-8pt]
$K\oplus k_v$, $K/k_v$ degree $2$ ram. field extn. & $()$ & $(12)$ & 2 & $1/2q$ \\[4pt]
\hline
 & & & &\\[-8pt]
$K\oplus k_v$, $K/k_v$ degree $2$ ram. field extn. & $(12)$ & $(12)$ & 2 & $1/2q$ \\[4pt]
\hline
 & & & &\\[-8pt]
$K$, $K/k_v$ degree $3$ ram. field extn. & $(123)$ & $(123)$ & 1 & $1/3q^2$ if $q\equiv 1 \pmod{3}$ \\[4pt]
 &  &  & & $0$ if $q\equiv 2 \pmod{3}$ \\[4pt]
\hline
 & & & &\\[-8pt]
$K$, $K/k_v$ degree $3$ ram. field extn. & $(132)$ & $(123)$ & 1 & $1/3q^2$ if $q\equiv 1 \pmod{3}$ \\[4pt]
 &  &  & & $0$ if $q\equiv 2 \pmod{3}$ \\[4pt]
\hline
 & & & &\\[-8pt]
$K$, $K/k_v$ degree $3$ ram. field extn. & $()$ & $(123)$ & 1 & $1/3q^2$ if $q\equiv 1 \pmod{3}$ \\[4pt]
 &  &  & & $0$ if $q\equiv 2 \pmod{3}$ \\[4pt]
\hline
 & & & &\\[-8pt]
$K$, $K/k_v$ degree $3$ ram. field extn. & $(12)$ & $(123)$ & 1 & $1/q^2$ if $q\equiv 2 \pmod{3}$ \\[4pt]
 &  &  & & $0$ if $q\equiv 1 \pmod{3}$ \\[4pt]
\hline
\end{tabular}

Let $S$ be the set of places of $k$ corresponding to the $\F_q$-rational points of $E$.  
Let $\Sigma$ be a choice $\Sigma_v$ of a cubic \'{e}tale $k_v$ algebra for each $v\in S$.
We define $c_\Sigma=\prod_{v\in S} c_{\Sigma_v}$, where the constants $c_{\Sigma_v}$ are defined by the chart above.
We also choose a $\Sigma'$ similarly.
Let $N_\Sigma(q^2n)$ be the number of isomorphism classes of cubic extensions $k'$ of $k$ such that for all $v\in S$, we have
$k' \tensor_k k_v \isom \Sigma_v$, and such that the norm of the relative discriminant $|D_{k'/k}|=q^{2n}.$
Now we can state the result of Datskovsky and Wright and cubic extensions of function fields.
\begin{theorem}[Corollary of Theorem 4.3 of \cite{DaWr}]\label{T:DW}
In the above notation
$$
\lim_{n\ra \infty} \frac{N_\Sigma(q^{2n})}{N_{\Sigma'}(q^{2n})}=\frac{c_\Sigma}{c_\Sigma'}.
$$ 
\end{theorem}

Theorem~\ref{T:EByPoint} now follows because the above chart gives all the possibilities for the completions of cubic extensions $k'$ of $k$ at each $\F_q$-rational place of $k$
(and the thus determined number of $\F_q$-rational points in that fiber), and Theorem~\ref{T:DW}
gives their relative probabilities.   For a triple cover $(C,\pi)\in T_{E,g}(\F_q)$ corresponding to an extension $k'/k$, we have  $|D_{k'/k}|=q^{2n}$, where
$g=n-1-\frac{3\chi(E)}{2}$, and thus we can replace the limit in $n$ with a limit in $g$.

\section{Further Directions and Conjectures}\label{S:Further}
The above computations of the distribution of points of trigonal curves suggests that similar questions could be attacked for $n$-fold covers of $\P^1$ (or an arbitrary curve) using methods from the study of counting extensions of global fields by their discriminants, which has been more heavily studied in the case of number fields.  
In analog to this work, one would naturally only study $n$-gonal curves with a specified monodromy group (i.e. the Galois group of the Galois closure of the field extension).
However, for this method to apply, one needs information on the densities of local behaviors of those number fields.  In the number field case, such results are only currently available for cubic extensions (the work of Datskovsky and Wright \cite{DaWr} used above, or going back to Davenport and Heilbronn \cite{S3} over $\Q$), for $(\Z/p \Z)^k$ Galois extensions (\cite{Wr, abelian}), and for degree $n$ extensions with Galois closure group $S_n$ and $n=4,5$) (by work of Bhargava \cite{S4,S5}).
The first two cases above also have results in the function field case, which have been applied in this paper and forthcoming work of the author, respectively.  It is intriguing to ask if a function field version of Bhargava's counting results \cite{S4,S5} could give the distribution of the number of points of $4$-gonal and $5$-gonal curves.  In analogy with the number field case, one expects
such methods would only count $4$-gonal curves with $S_4$ monodromy, and that $D_4$-monodromy curves would be a non-negligible portion of all $4$-gonal curves.  Since the methods for counting $D_4$ quartic extensions have not given any results on densities of local behaviors or independence of those behaviors, understanding the distribution of points of $4$-gonal curves with $D_4$ monodromy would be very interesting.

	In the case of $n$-gonal curves with full $S_n$-monodromy, we can at least predict what the distribution of the number of points in each fiber should be using the heuristics of Bhargava \cite[Conjecture 5.1]{BhH}.  In particular, we give the prediction for the average.
\begin{conjecture}\label{C}
 Let $E$ be a smooth curve over a finite field $\F_q$ of characteristic $> n$ and fix $z\in E(\F_q)$.
The average number of points in the fiber over $z$ of a 
random curve over $\F_q$ with a degree $n$ map to $E$ and full monodromy (in the $g\ra \infty$ limit) is 
$$
\frac{\sum_k k \sum_{\ell=0}^{n-1} \frac{p(n,n-\ell,k)}{q^\ell}}{\sum_{\ell=0}^{n-1} \frac{p(n,n-\ell)}{q^\ell}}=
1+\frac{1}{q}+O(\frac{1}{q^2}),
$$
where $p(n,m,k)$ is the number of partitions of $n$ into $m$ parts such that the parts take exactly $k$ values and $p(n,m)$ is the number of partitions of $n$ into $m$ parts.
\end{conjecture}
It seems within reach to prove this conjecture for $n=4,5$ by proving function field analogs of \cite{S4,S5}.
For $n=4$ the expected fiber size is (conjecturally) $$
1+\frac{q^2+q}{q^3+q^2+2q+1}
$$
and for $n=5$ the expected fiber size is (conjecturally)
$$
1+\frac{q^3+2q^2+2q}{q^4+q^3+2q^2+2q+1}.
$$
 In particular, we now show that Conjecture~\ref{C}
follows from the function field analog of \cite[Conjecture 5.1]{BhH}.

Let $k$ be the function field of $E$ and $v$ place of $E$ corresponding to a $k_v$-rational point.  Then \cite[Conjecture 5.1]{BhH}.gives conjectures for local densties for $k_v$ algebras
among $k' \tensor_k k_v$ while $k'$ ranges over degree $n$ $S_n$-extensions of $k$.  In particular, \cite[Conjecture 5.1]{BhH} predicts, as we range over degree $n$ \'etale $k_v$ algebras $L$, that $L$ 
appears with relative density $\frac{1}{|\Aut(L)||D_{L/{k_v}}|}$ (recall $|D_{L/{k_v}}|$ is the absolute norm of the relative discriminant of $L/{k_v}$).
Recall that isomorphism classes of degree $n$ \'etale $k_v$ algebras exactly correspond to continuous homomorphisms $\Gal(\bar{k_v}/k_v)\ra S_n$, and we can rephrase the above to say that
$L$  appears with relative density 
$$\frac{\#|\chi: \Gal(\bar{k_v}/k_v)\ra S_n \textrm{ corr. to } L|}{|D_{L/{k_v}}|}.$$
Continuous homomorphisms $\Gal(\bar{k_v}/k_v)\ra S_n$ correspond to choices $x,y\in S_n$ such that $xyx^{-1}=y^q$, and $|D_{L/{k_v}}|=q^{n-\#\operatorname{cycles}(y)}$.
The number of $\F_q$ rational points above $v$ in the curve corresponding to the extension $k'$ of $k$ is given by the number of 
$\langle x,y \rangle$-orbits on $\{1,\dots,n\}$ that are also $\langle y \rangle$-orbits.

Let $y\in S_n$ have $a_i$ orbits of size $b_i$, for $i=1,\dots k$.  Then the number of $x\in S_n$ such that $xyx^{-1}=y^q$ is $\prod_ia_i! b_i^{a_i}$.  
If we consider a single cycle $\sigma$ of $y$ of length $b_i$, then the number of $x\in S_n$ such that $\sigma$ remains a $\langle x,y \rangle$-orbit is
$\frac{1}{a_i}\prod_i a_i! b_i^{a_i}$.  So given the choice of $y$, the expected contribution of $\sigma$ to rational points is $\frac{1}{a_i}$.  Thus, given a choice of $y$,
the expected number of rational points is $k$, the number of distinct cycle lengths of $y$.
Note that $\prod_ia_i! b_i^{a_i}$ is the size of the centralizer of $y$, so choosing a permutation $y\in S_n$ with relative probability  $\frac{\prod_ia_i! b_i^{a_i}}{q^{n-\#\operatorname{cycles}(y)}}$
is the same as choosing it with relative probability $\frac{|\operatorname{Cent}(y)|}{q^{n-\#\operatorname{cycles}(y)}}$ or $\frac{1}{|\operatorname{ConjClass}(y)|q^{n-\#\operatorname{cycles}(y)}}$.
This is equivalent to choosing a conjugacy class in $S_n$, i.e. a partition of $n$, with relative probability $\frac{1}{q^{n-\#\operatorname{parts}}}$.
Since a partition with $k$ distinct part sizes, gives an expected $k$ parts, Conjecture~\ref{C} follows from the function field analog of \cite[Conjecture 5.1]{BhH}.


\begin{thebibliography}{99}

\bibitem{S4} M. Bhargava, The density of discriminants of quartic rings and fields, Ann. of Math. (2) {\bf 162} (2005), 1031--1063. 

\bibitem{S5} M. Bhargava, The density of discriminants of quintic rings and fields,  Ann. of Math., (3) {\bf 172} (2010),  1559--1591.

\bibitem{BhH} M. Bhargava, Mass formulae for extensions of local fields, and conjectures on the density of number field discriminants, Int. Math. Res. Not. IMRN {\bf 2007}, no.~17, Art. ID rnm052, 20 pp.

\bibitem{BuDaFeLa1}  A.~Bucur, C.~David, B.~Feigon and M.~Lal\'in, Statistics for traces of cyclic trigonal curves over finite fields,  \textit{Int. Math. Res. Not. IMRN},  Advance Access published on October 27, 2009, doi:10.1093/imrn/rnp162.

\bibitem{BuDaFeLa2} A.~Bucur, C.~David, B.~Feigon and M.~Lal\'in, Biased statistics for traces of cyclic p-fold covers over finite fields, to appear in Proceedings of 2008 Banff Workshop "Women In Numbers'' (in the Fields Institute Communications Series).

\bibitem{BuDaFeLa3} A.~Bucur, C.~David, B.~Feigon and M.~Lal\'in, Fluctuations in the number of points on smooth plane curves over finite fields, 
\textit{Journal of Number Theory} \textbf{130} (2010), 2528--2541.

\bibitem{BuKe} A.~Bucur, and K.~Kedlaya, 
The probability that a complete intersection is smooth, arXiv:1003.5222v1

\bibitem{DaWr} B. Datskovsky\ and\ D. J. Wright, The adelic zeta function associated to the space of binary cubic forms. II. Local theory, J. Reine Angew. Math. {\bf 367} (1986), 27--75.

\bibitem{S3}
H. Davenport\ and\ H. Heilbronn, On the density of discriminants of cubic fields. II, Proc. Roy. Soc. London Ser. A {\bf 322} (1971), no.~1551, 405--420.


\bibitem{ks} N.~M.~Katz and P.~Sarnak, \textit{Random matrices, Frobenius eigenvalues, and monodromy}. American Mathematical Society Colloquium Publications, 45. American Mathematical Society, Providence, RI, 1999. xii+419 pp.

\bibitem{KuRu} P.~Kurlberg and Z.~Rudnick,  The fluctuations in the number of points on a hyperelliptic curve over a finite field, \textit{J. Number Theory} \textbf{129}  (2009),  no. 3, 580--587.

\bibitem{KuWi} P.~Kurlberg and I.~Wigman. Gaussian point count statistics for families of curves over a fixed finite field. \textit{Int. Math. Res. Not. IMRN}, to appear.

\bibitem{abelian} M.~M.~Wood, On the probabilities of local behaviors in abelian field extensions, \textit{Compos. Math.} \textbf{146} (2010), no. 1, 102--128.

\bibitem{Wr}
D. J. Wright, Distribution of discriminants of abelian extensions, Proc. London Math. Soc. (3) {\bf 58} (1989), no.~1, 17--50. 


\end{thebibliography}
\end{document}